\documentclass[%
 print,
superscriptaddress,
 amsmath,amssymb,
 aps,
]{revtex4-1}

\usepackage{graphicx}
\usepackage{dcolumn}
\usepackage{bm}
\usepackage[french]{babel}

\begin{document}


\title{A new mathematical symbol : the termirial }

\author{Claude-Alexandre Simonetti}
\affiliation{LPC Caen}
\affiliation{EAMEA Cherbourg}
\affiliation{ENSICAEN}
\date{\today}

\def\termi{\overset{|}{\raisebox{-1.0ex}{\tiny{+}}}}
\def\termisecond{\overset{(2)}{\termi}}
\def\termiminusone{\overset{(-1)}{\termi}}
\def\termizero{\overset{(0)}{\termi}}
\def\termip{\overset{(p)}{\termi}}
\def\termii{\overset{(i)}{\termi}}
\def\termipi1{\overset{(p-i-1)}{\termi}}
\def\termipminusi{\overset{(p+1-i-1)}{\termi}}
\def\termiiplusone{\overset{(i+1)}{\termi}}
\def\termipplusone{\overset{(p+1)}{\termi}}
\def\termipminusone{\overset{(p-1)}{\termi}}

\begin{abstract}

The understanding of probability can be difficult for a few young scientists. Consequently, this new mathematical symbol, related to binomial coefficients and simplicial polytopic numbers, could be helpful to science education. Moreover, one can obtain kinds of remarkable identities and generalize them to a sort of ``Newton's binomial theorem''. Finally, this symbol could be perhaps useful to other scientific subjects as well, such as computer science.

\end{abstract}

\pacs{Valid PACS appear here}
\maketitle


\section{\label{sec:level1}Study of binomial coefficients}
\subsection{Reminder}
In a set of n elements, the number of combinations of parts of p elements is the following:\\
\begin{equation} \label{coef_binom1}
\binom{n}{p}=\frac{n!}{p!~(n-p)!}
\end{equation}\\
With $\binom{n}{p}$ the binomial coefficient, $(n,p)\in{\mathbb{N}^2}$ and $p\le n$. And the exclamation point in equation \eqref{coef_binom1} is named ``factorial'' and is defined as, $\forall n\in{\mathbb{N}^*}$:\\
\begin{eqnarray} \nonumber
n!&=&n\cdot (n-1) \cdot (n-2)~...~3 \cdot 2 \cdot 1\\ 
n!&=&\prod_{i=1}^{n} i \label{factorielle}
\end{eqnarray}
For $n=0$, by definition : 0!=1.

\subsection{From factorial to termirial}
For its part, the termirial is defined as, $\forall n\in{\mathbb{N^*}}$:\\
\begin{eqnarray} \nonumber
n\termi&=&n + (n-1) +...+3 + 2 + 1\\ 
n\termi&=&\sum_{i=1}^{n} i\label{termirielle1}
\end{eqnarray}\\
This is called a triangular number of order n, see reference \cite{9780486442327}. The termirial symbol is almost like the factorial symbol, with just a little difference, however : the dot ($\cdot$) of the interrogation point, which could be a reminder of a multiplication, is replaced by a ``plus'' ($+$) sign. Indeed, instead of multiplying factors of a multiplication -- factorial -- one adds terms of an addition, hence the name of ``termirial''. Up to this point, things look quite basic. However, as the termirial is the $n^{th}$ partial sum of an arithmetic sequence $(U_n)_{n\in{\mathbb{N^*}}}$ with $U_1=1$ as first term and $r=1$ as common difference, one can notice a first thing, $\forall n\in{\mathbb{N^*}}$:\\
\begin{eqnarray} \label{termirielle2}
n\termi&=&\frac{n\cdot (n+1)}{2}=\binom{n+1}{2}=\binom{n+1}{n-1}
\end{eqnarray}\\
One will get back to it later, with the ``generalized'' termirial.

\subsection{Remarkable identity}
It is possible to obtain a kind of remarkable identity. Indeed, $\forall (n,m)\in{\mathbb{N^*}^2}$:
\begin{eqnarray} \nonumber
(n+m)\termi&=&\frac{(n+m)\cdot (n+m+1)}{2}\\ \nonumber
(n+m)\termi&=&\frac{n\cdot (n+1)}{2}+\frac{2\cdot n \cdot m}{2}+\frac{m\cdot (m+1)}{2}\\ \label{idrem1}
(n+m)\termi&=&n\termi +n\cdot m+ m\termi 
\end{eqnarray}

\subsection{Intellectual path}
\begin{figure}[!h]
\begin{center}
\includegraphics*[scale=0.6]{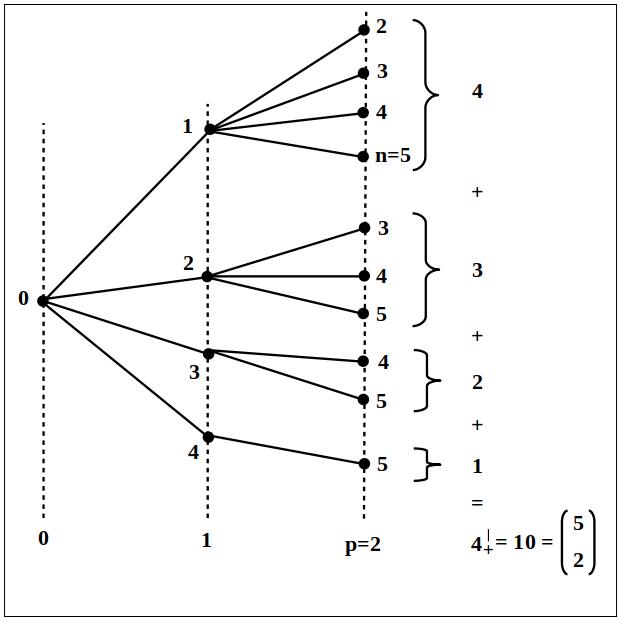}
\caption{This figure represents the intellectual path from the binomial coefficient $\binom{5}{2}$ to the termirial of 4, through a classical tree view. Here, only the energy level $m_1$ is represented, which is sufficient, because one only has 2 energy levels, not more.}
\label{fig:Termi_1}
\end{center}
\end{figure}

At the beginning, I was studying the black body radiation (Max Planck's law). In this context, suppose one has $n$ particles to fill up in 2 discrete energy levels determined by quantum mechanics. Suppose that the first energy level named $m_1$ can contain $p$ particles, and that the second one ($m_2$) can contain $(n-p)$ particles. The number of possible combinations is same as equation \eqref{coef_binom1}:\\
\begin{equation} \label{coef_binom2}
\binom{n}{p}=\frac{n!}{p!~(n-p)!}
\end{equation}\\
Let us start with a simple example : $n=5$ and $p=2$. The figure \ref{fig:Termi_1} tries to explain the intellectual path which has led to the termirial. One can see the energy level $m_1$ containing $p=2$ particles, with all the ways to fill it up, for particles particules ranked from 1 to 5. In this case, as a reminder, it is an unordered sampling without replacement, the rank order of the particles has no importance : $\{1~2\}$ and $\{2~1\}$ are counted only once. Finally :  $\binom{5}{2}=4 \termi = 10$.

\begin{figure}[!h]
\begin{center}
\includegraphics*[scale=0.6]{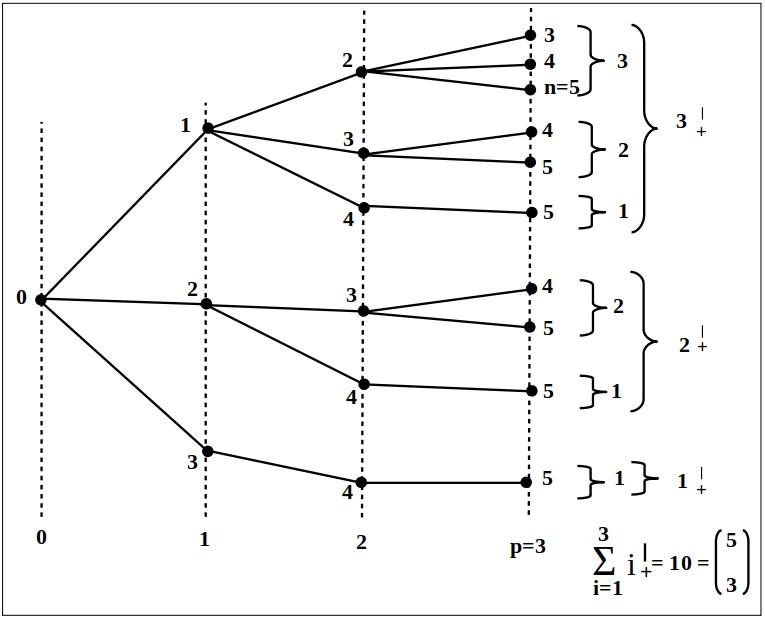}
\caption{This figure represents the intellectual path from the binomial coefficient $\binom{5}{3}$ to a sum of termirials, a kind of ``termirial of termirials'' or a $2^{nd}$ termirial. From now, one can generalize to the $3^{rd}$ termirial, the $4^{th}$ termirial, etc.}
\label{fig:Termi_2}
\end{center}
\end{figure}

Let us make things a little bit more complicated, with $n=5$ and $p=3$. In terms of binomial coefficients, one has the same result as before : $\binom{5}{3}=\binom{5}{2}=10$. But the intellectual path is a bit different, see figure \ref{fig:Termi_2}.

One can see a kind of ``termirial of termirials'', which is called a tetrahedral number of order n, see reference \cite{9780486442327}:\\
\begin{eqnarray} \nonumber
n\termisecond&=&n\termi + (n-1)\termi +...+3\termi+ 2\termi + 1\termi\\ 
n\termisecond&=&\sum_{k=1}^{n} k \termi =\sum_{k=1}^{n} \sum_{i=1}^{k} i=\sum_{k=1}^{n}\frac{k(k+1)}{2!}\label{termirielle3}\\
n\termisecond&=&\binom{n+2}{3}=\binom{n+2}{n-1}=\frac{n(n+1)(n+2)}{3!}\label{termirielle4}
\end{eqnarray}\\
The transition from \eqref{termirielle3} to \eqref{termirielle4} can be done with mathematical induction. Indeed, if the following proposition $P(n)$ is true:
\begin{eqnarray} \label{termirielle5}
n\termisecond&=&\sum_{k=1}^{n}\frac{k(k+1)}{2!}=\frac{n(n+1)(n+2)}{3!}
\end{eqnarray}
And, for $P(1)$:
\begin{eqnarray} \nonumber
1\termisecond&=&\sum_{k=1}^{1}\frac{k(k+1)}{2}=\frac{1 \cdot (1+1)}{2}=1\\
\binom{3}{3}&=&\frac{1(1+1)(1+2)}{3!}=1\nonumber
\end{eqnarray}
Finally, with $P(n+1)$:
\begin{eqnarray} \nonumber
(n+1)\termisecond&=&\sum_{k=1}^{n}\frac{k(k+1)}{2}+\frac{(n+1)(n+2)}{2}\\
                 &=&\frac{n(n+1)(n+2)}{3!}+\frac{(n+1)(n+2)}{2}\nonumber\\
                 &=&\frac{n(n+1)(n+2)}{6}+3\frac{(n+1)(n+2)}{6}\nonumber\\
(n+1)\termisecond&=&\frac{(n+1)(n+2)(n+3)}{3!}\label{termirielle6}
\end{eqnarray}
\textit{QED}.
\subsection{Remarkable identity of the $2^{nd}$ termirial}
Like the $1^{st}$ termirial, it is possible to obtain a kind of remarkable identity. Indeed, $\forall (n,m)\in{\mathbb{N^*}^2}$:
\begin{eqnarray} \nonumber
(n+m)\termisecond&=&\frac{(n+m+2)}{3} \cdot (n+m)\termi\\ \nonumber
                 &=&\frac{(n+m+2)}{3} \cdot (n\termi +n\cdot m+ m\termi) \\ \label{idrem2}
(n+m)\termisecond&=&n\termisecond +n\cdot m\termi +m\cdot n\termi + m\termisecond 
\end{eqnarray}\\
From now on, one can generalized to the $p^{th}$ termirial.

\section{\label{sec:level2}Generalization of the termirial}
\subsection{The $p^{th}$ termirial}
Obviously, the the $p^{th}$ termirial, which is called a (p+1)-simplicial polytopic number, see reference \cite{9780486442327}, will be defined as, $\forall (n,p)\in{\mathbb{N^*}\times \mathbb{N^*}}$:\\
\begin{eqnarray} \label{termirielle7}
n\termip&=&\sum_{m=1}^{n} \sum_{l=1}^{m}...\sum_{j=1}^{k} \sum_{i=1}^{j} i\\
n\termip&=&\sum_{k=1}^{n} k\termipminusone=\sum_{k=1}^{n}\frac{1}{p!}\prod_{i=0}^{p-1} (k+i)\label{termirielle8}\\
n\termip&=&\binom{n+p}{p+1}=\binom{n+p}{n-1}\label{termirielle9}\\
n\termip&=&\frac{1}{(p+1)!} \prod_{i=0}^{p} (n+i) \label{termirielle10}
\end{eqnarray}\\
With, for the equation \eqref{termirielle7}, a number of \textit{sigmas} $(\sum)$ symbols  equal to $p$.\\
For $p=0$, one can define the $0^{th}$ termirial as:
\begin{equation} \nonumber
n\termizero=n=\binom{n}{1}
\end{equation}
For $p=-1$, one can define the $-1^{th}$ termirial as:
\begin{equation} \nonumber
n\termiminusone=1=\binom{n-1}{0}
\end{equation}
The proof to go from \eqref{termirielle8} to \eqref{termirielle9} (or from \eqref{termirielle8} to \eqref{termirielle10}) can be done with mathematical induction as well. Indeed, admitting that the following $P(n)$ proposition is true:
\begin{eqnarray} \label{termirielle11}
\sum_{k=1}^{n}\frac{1}{p!}\prod_{i=0}^{p} (k+i)&=&\frac{1}{(p+1)!} \prod_{i=0}^{p} (n+i)
\end{eqnarray}
One have, for $P(1)$:
\begin{eqnarray} \nonumber
1\termip=\frac{1}{p!}\prod_{i=0}^{p-1} (1+i)=\frac{p!}{p!}&=&1\\
\frac{1}{(p+1)!} \prod_{i=0}^{p} (1+i)=\frac{(p+1)!}{(p+1)!}&=&1
\end{eqnarray}
And, for $P(n+1)$:
\begin{eqnarray} \nonumber
(n+1)\termip&=&\sum_{k=1}^{n}\frac{1}{p!}\prod_{i=0}^{p-1} (k+i)+\frac{1}{p!}\prod_{i=0}^{p-1} (n+1+i)\\
                 &=&\frac{1}{(p+1)!} \prod_{i=0}^{p} (n+i)+\frac{(p+1)}{(p+1)!}\prod_{i=1}^{p} (n+i)\nonumber\\
                 &=&\frac{(n+p+1)}{(p+1)!} \prod_{i=1}^{p} (n+i)\nonumber\\
(n+1)\termip&=&\frac{1}{(p+1)!} \prod_{i=1}^{p+1} (n+i)=\binom{n+p+1}{p+1}\label{termirielle12}
\end{eqnarray}
\textit{QED}.

\subsection{Pascal's rule}
With Pascal's rule as a reminder:
\begin{eqnarray} \nonumber
\binom{n}{p}+\binom{n}{p+1}&=&\binom{n+1}{p+1}
\end{eqnarray}\\
It is possible to adapt it to termirials:
\begin{eqnarray} \nonumber
(n+1)\termip+(n)\termipplusone&=&(n+1)\termipplusone
\end{eqnarray}

\subsection{Newton's binomial theorem}
With Newton's binomial theorem as a reminder, $\forall (a,b)\in{\mathbb{R}\times \mathbb{R}}$ and $\forall n\in{\mathbb{N}}$:
\begin{eqnarray} \nonumber
(a+b)^n&=&\sum_{p=0}^{n} \binom{n}{p} a^{n-p} \cdot b^p 
\end{eqnarray}
It is possible to adapt it as well to termirials, $\forall (n,m)\in{\mathbb{N}^* \times \mathbb{N}^*}$ and $\forall p\in{\{\{-1\} \cup \mathbb{N}\}}$:
\begin{eqnarray} \label{Newton}
(n+m)\termip&=&\sum_{i=-1}^{p} n\termii \cdot m\termipi1
\end{eqnarray}
The proof can be done with mathematical induction. Indeed, admitting that the following $P(p)$ proposition is true:
\begin{eqnarray} \nonumber
(n+m)\termip&=&\sum_{i=-1}^{p} n\termii \cdot m\termipi1
\end{eqnarray}
First of all, equations \eqref{idrem1} and \eqref{idrem2} respectively hold $P(1)$ and $P(2)$, even if only $P(1)$ -- or only $P(2)$ -- is enough for this mathematical induction. Considering $P(p+1)$:
\begin{eqnarray} \nonumber
(n+m)\termipplusone&=&\frac{n+m+p+1}{p+2} (n+m)\termip \\ \nonumber
                   &=&\frac{n+m+p+1}{p+2} \sum_{i=-1}^{p} n\termii \cdot m\termipi1 \\ \nonumber
                   &=&\frac{n+m+p+1}{p+2} \sum_{i=-1}^{p} n\termii \cdot m\termipi1 \\ \nonumber
                   &=&\frac{1}{p+2} \sum_{i=-1}^{p} \{[n+(i+1)]+[m+p+1-(i+1)]\}\cdot n\termii \cdot m\termipi1 \\ \nonumber
                   &=&\frac{1}{p+2} \sum_{i=-1}^{p} [n+(i+1)]\cdot n\termii \cdot m\termipi1 +\frac{1}{p+2} \sum_{i=-1}^{p} [m+p+1-(i+1)]\cdot n\termii \cdot m\termipi1 \\ \nonumber
                   &=&\frac{1}{p+2} \sum_{i=-1}^{p} (i+2)\cdot n\termiiplusone \cdot m\termipi1 +\frac{1}{p+2} \sum_{i=-1}^{p} [p+2-(i+1)]\cdot n\termii \cdot m\termipminusi \\ \nonumber
                   &=&\frac{1}{p+2} \sum_{i=0}^{p+1} (i+1) n\termii m\termipminusi + \sum_{i=-1}^{p} n\termii m\termipminusi -\frac{1}{p+2} \sum_{i=-1}^{p} (i+1) n\termii m\termipminusi \\ \nonumber
                   &=&n\termipplusone \cdot m\termiminusone + \sum_{i=-1}^{p} n\termii m\termipminusi \\ \nonumber
(n+m)\termipplusone&=&\sum_{i=-1}^{p+1} n\termii m\termipminusi \\ \nonumber
\end{eqnarray}
\textit{QED}.

\section{\label{sec:level3}Pratical applications}
\subsection{Computational complexity}
In computer science, the termirial can be used to calculate the complexity of computer programs with intricated ``for'' loops as follows : for example, with 4 ``for'' loops, with:\\
- for i from 1 to n=100~;\\
- for j from 1 to i~;\\
- for k from 1 to j~;\\
- for l from 1 to k.\\
\\
The complexity will be:
\begin{eqnarray} \nonumber
100\overset{(4-1)}{\termi}&=&\binom{100+3}{3+1}=\binom{103}{4}=\binom{103}{99}=4~421~275
\end{eqnarray}
So, approximatey 4,4 millions of operations.\\
More generally, the complexity of this kind of program will be, according to equation \eqref{termirielle10} and $\forall (n,p)\in {\mathbb{N^*}\times \mathbb{N^*}}$ : $n\termipminusone=\Theta (n^{p})$.
\subsection{Pseudo-fractal aspect}

\begin{figure}[!h]
\begin{center}
\includegraphics*[scale=0.6]{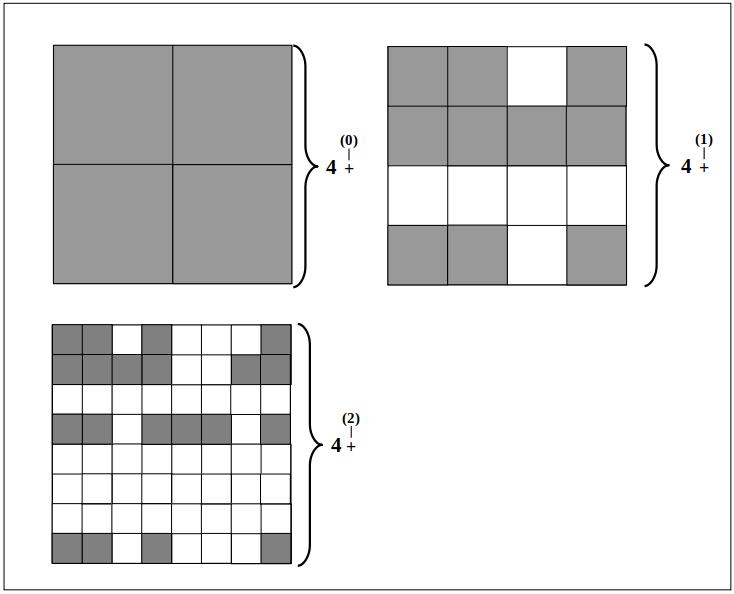}
\caption{This figure tries to convert $p^{th}$ termirial of 4 into fractal figures, for p=0, p=1 and p=2. Only the grey squares are counted.}
\label{fig:Termi_3}
\end{center}
\end{figure}

The figure \ref{fig:Termi_3} tries to convert $p^{th}$ termirial of 4 into fractal figures, for p=0, p=1 and p=2. Only the grey squares are counted. For example, for $4\termisecond$, 20 squares are in grey, which corresponds to $4\termisecond=20$. For each incrementation of p, the side of each square is divided by 2, and the same basic pattern is repeated:\\
- $4\termi$ when the 4 previous squares are in grey~;\\ 
- $3\termi$ when the 3 previous squares are in grey~;\\
- $2\termi$ when the 2 previous squares are in grey~;\\
- $1\termi$ when only 1 previous square is in grey.\\

But the Hausdorff dimension (D) of this presumed fractal is not constant, and the limit of D as p approaches infinity is 2. Indeed, let $a$ the side length of each square related to the $p^{th}$ termirial. The total surface $S_p$ of the squares in grey will be:\\

\begin{eqnarray} \nonumber
S_p&=&a^2 \cdot n\termip = a^2 \cdot \frac{1}{(p+1)!} \prod_{i=0}^{p} (n+i)
\end{eqnarray}
And for the $S_{p-1}$ surface, related to the $(p-1)^{th}$ termirial, for a side length of each square doubled in size:\\
\begin{eqnarray} \nonumber
S_{p-1}&=&4a^2 \cdot n\termipminusone = 4a^2 \cdot \frac{1}{(1)!} \prod_{i=0}^{p-1} (n+i)
\end{eqnarray}
Hence, for $p > 0$:\\
\begin{eqnarray} \nonumber
\frac{S_{p-1}}{S_{p}}&=&4 \cdot \frac{p+n}{p+1}= 4 \cdot \frac{1+n/p}{1+1/p}
\end{eqnarray}
Finally, this surface ratio, as p approaches infinity, is 4:\\
\begin{eqnarray} \nonumber
\lim_{p \to \infty} \frac{S_{p-1}}{S_{p}}&=&4
\end{eqnarray}
Consequently, the Hausdorff dimension (D) is 2, as $2^D=4$, when p approaches infinity. To conclude, as D is not constant, even if figure \ref{fig:Termi_3} looks interesting, termirials do not seem to be related to fractals.

\section{Conclusion}
The termirial is a symbol which could be helpful to scientific students, for the understanding of probability. As well, it could be perhaps helpful in other subjects than science education, that are out of this pr\'e-publication.

\nocite{*}

\bibliography{Termirial_v4.bib}

\end{document}